\documentclass[runningheads]{llncs}
\usepackage{amsmath}
\usepackage{amssymb}
\usepackage[applemac]{inputenc}
\usepackage[colorlinks]{hyperref}
\usepackage{url}
\usepackage{graphicx,epstopdf}
\usepackage{subcaption}
\usepackage[T1]{fontenc}
\newcommand{\normal}{\vec{n}}
\newcommand{\tangent}{\vec{t}}
\newcommand{\Nor}{\operatorname{Nor}}

\newcommand{\Ver}{\operatorname{Ver}}
\usepackage{graphicx}
\begin{document}
\title{3 methods to put a Riemannian metric on Shape~Space\thanks{Supported by by FWF grant I 5015-N, Institut CNRS Pauli and University of Lille}}
%
%

\author{Alice~Barbora~Tumpach\inst{1} \and
Stephen~C.~Preston\inst{2}}
\authorrunning{Tumpach and Preston}
%
\institute{
Institut CNRS Pauli, Vienna, Austria and
University of Lille, France
\email{alice-barbora.tumpach@univ-lille.fr}\\
\and
Brooklyn College and CUNY Graduate Center, New York USA \\
\email{stephen.preston@brooklyn.cuny.edu}\\}

\maketitle              
\begin{abstract}
In many applications, one is interested in the shape of an object, like the contour of a bone or the trajectory of joints of a tennis player, irrespective of the way these shapes are parameterized. However for analysis of these shape spaces, it is sometimes useful to have a parameterization at hand, in particular if one is interested in deforming shapes. The purpose of the paper is to examine three different methods that one can follow to endow shape spaces with a Riemannian metric that is measuring deformations in a parameterization independent way.

\end{abstract}
 \vspace{-1cm}



\begin{figure*}[!ht]
 		\centering
\begin{tabular}{|p{2cm}|p{3.5cm}|p{3.5cm}|p{3.5cm}|}
\hline
Group $G$& Elements of an orbit under $G$ & one preferred element in the orbit & another choice of preferred element \\
\hline
\begin{center} $\mathbb{R}^2$ acting by translations \end{center}& \begin{center}\includegraphics[width=3cm]{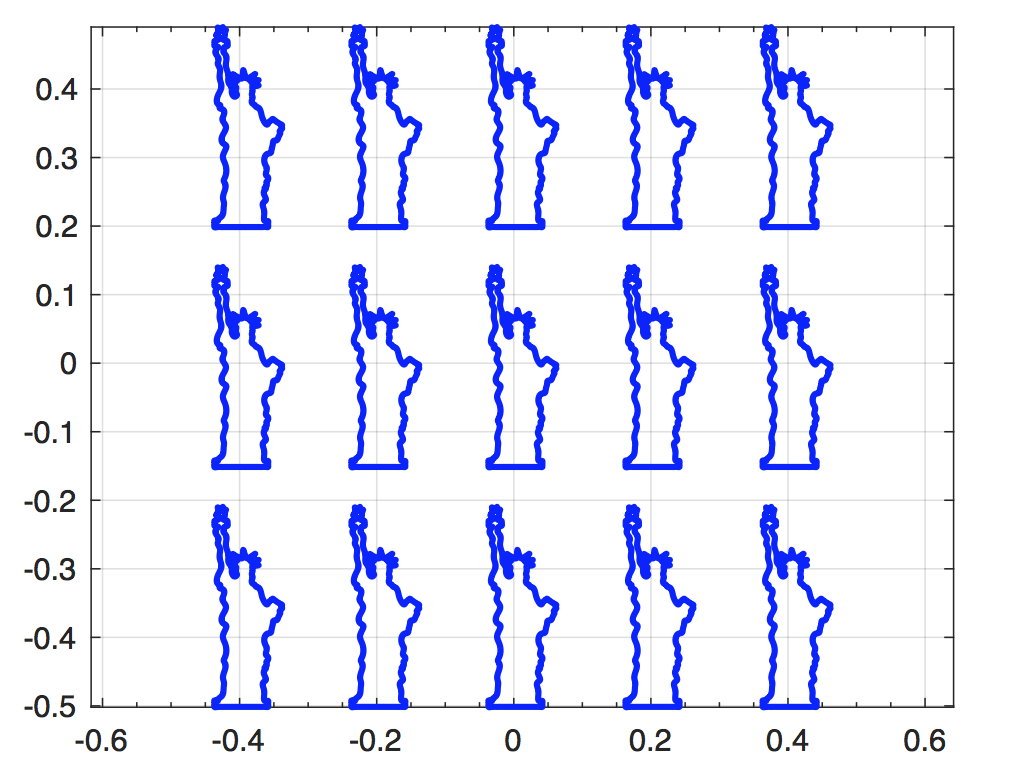}\end{center} & \begin{center}\includegraphics[width=2.6cm]{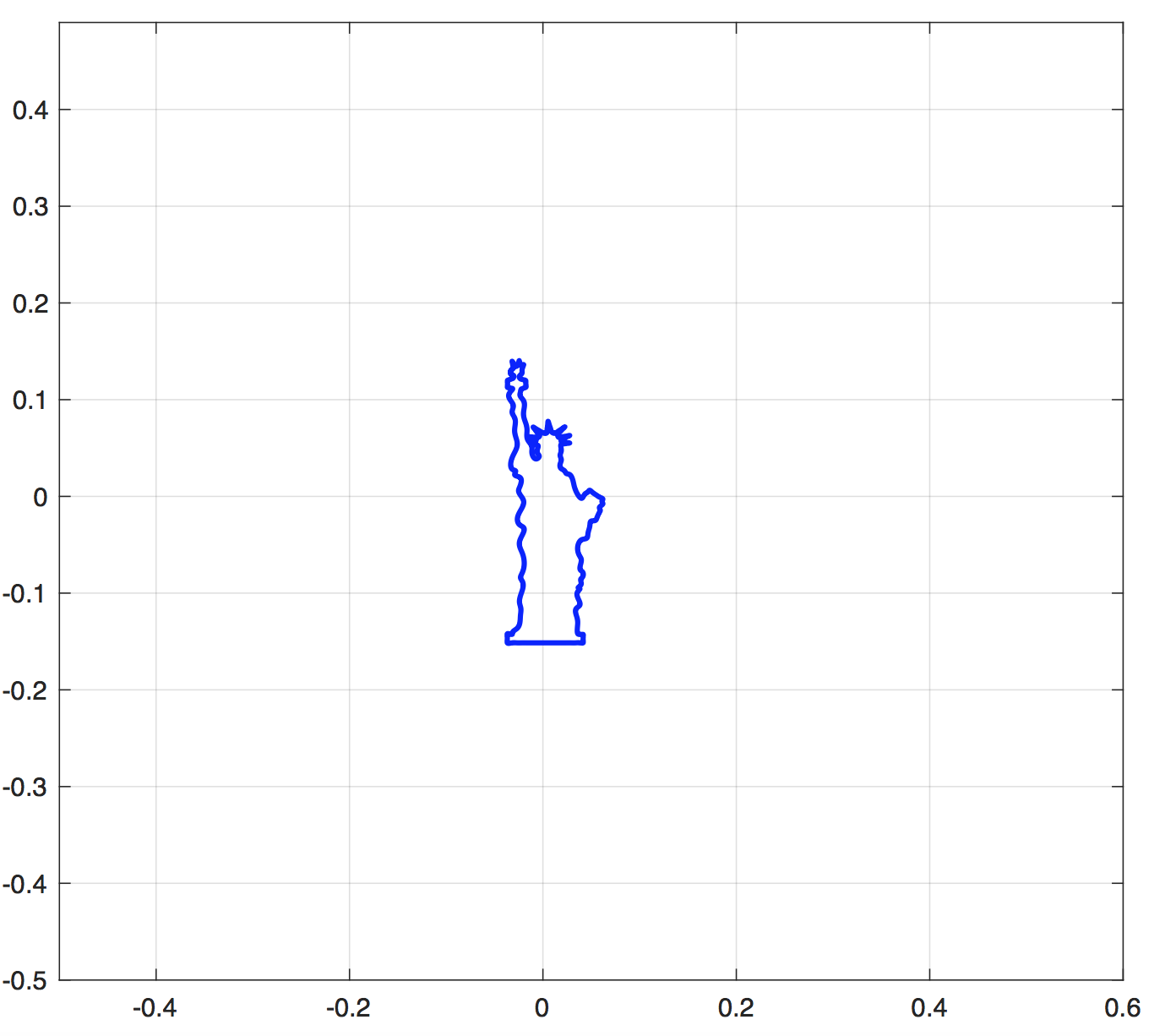} {\small \qquad centroid at origin}\end{center} 
& \begin{center}\includegraphics[width=2.8cm]{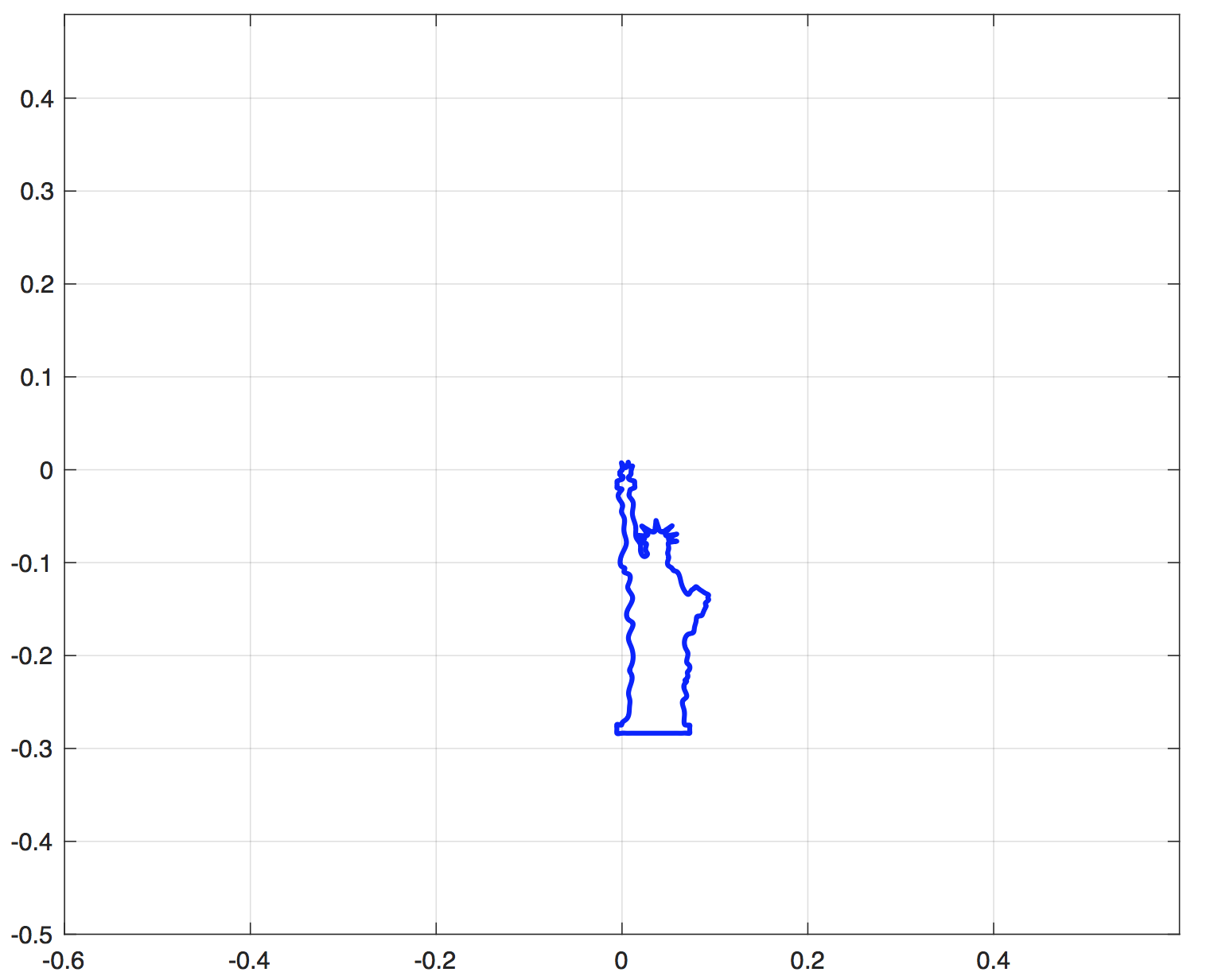} {\small curve starting at origin}\end{center}\\ 
\hline
\begin{center} $\mathrm{SO}(2)$ acting by rotations \end{center} &\begin{center}\includegraphics[width=3.2cm]{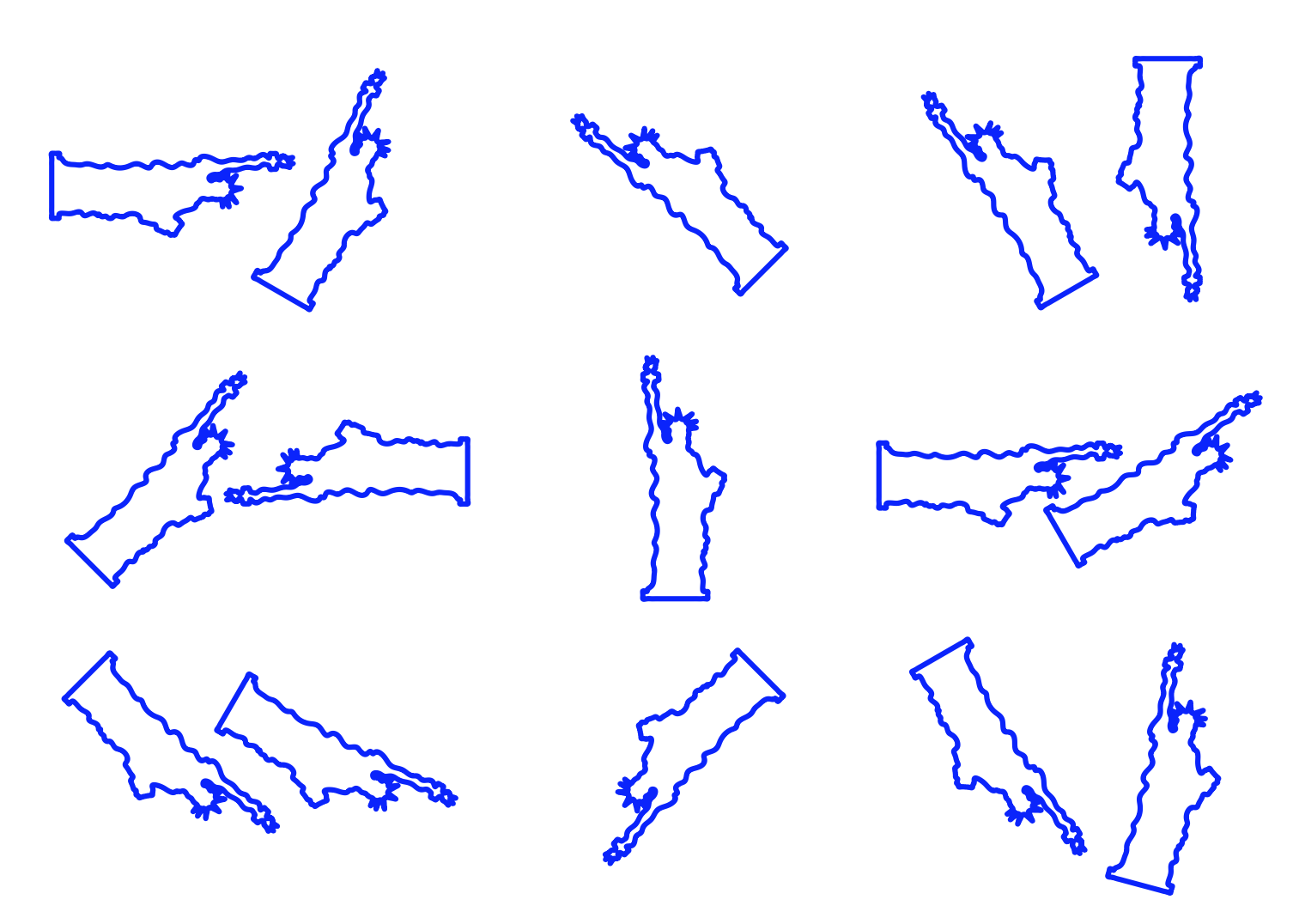}\end{center} & \begin{center} \begin{tabular}{l} \includegraphics[height=0.8cm]{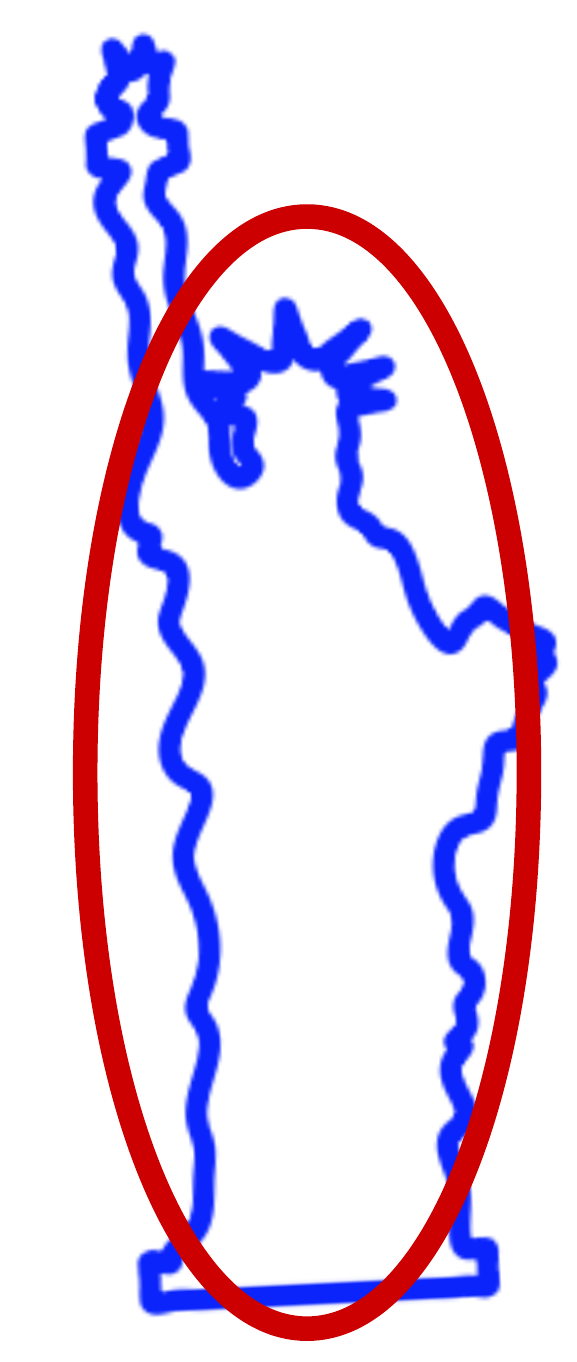}\\ axes of \\ approximating\\ ellipse aligned \end{tabular}\end{center}& \begin{center} \begin{tabular}{l} \includegraphics[width=0.8cm]{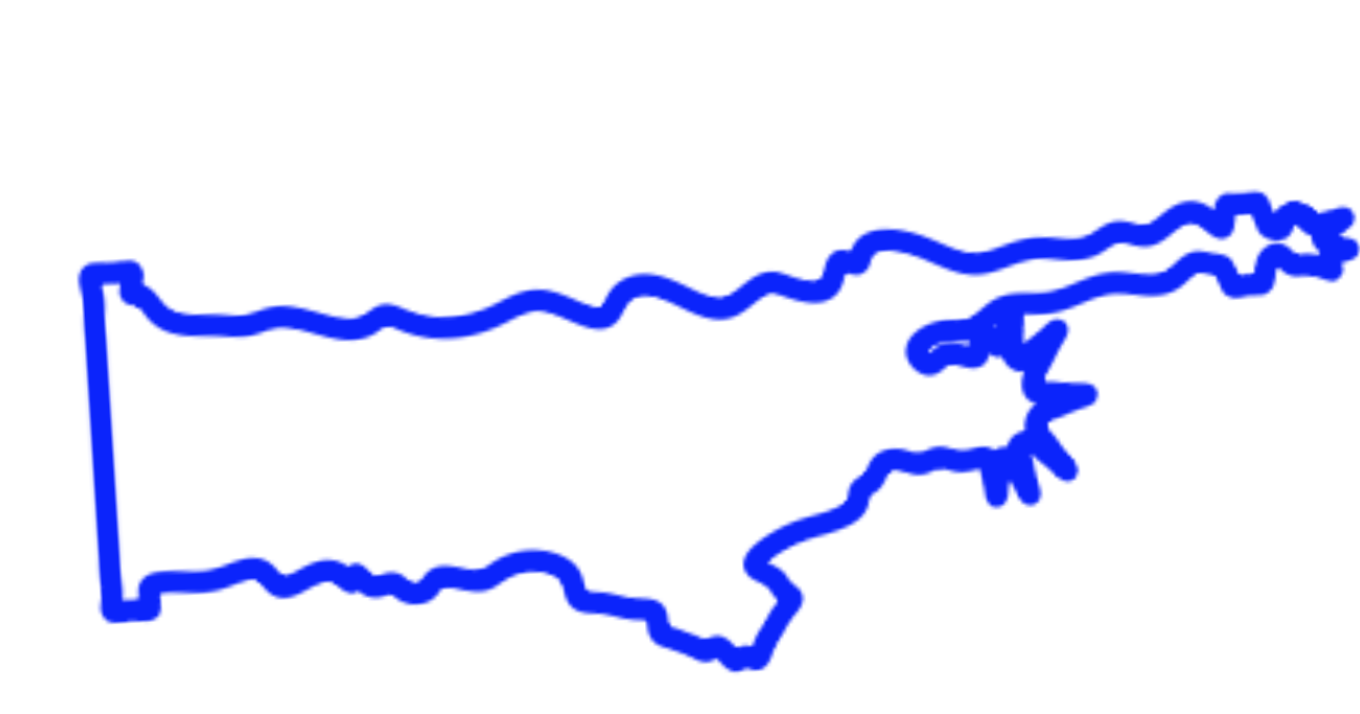} \\   tangent vector\\ at starting point\\ horizontal \end{tabular}\end{center} \\
\hline
%
\begin{center} $\mathbb{R}^+$ acting by scaling \end{center}&\begin{center} \includegraphics[width=3cm]{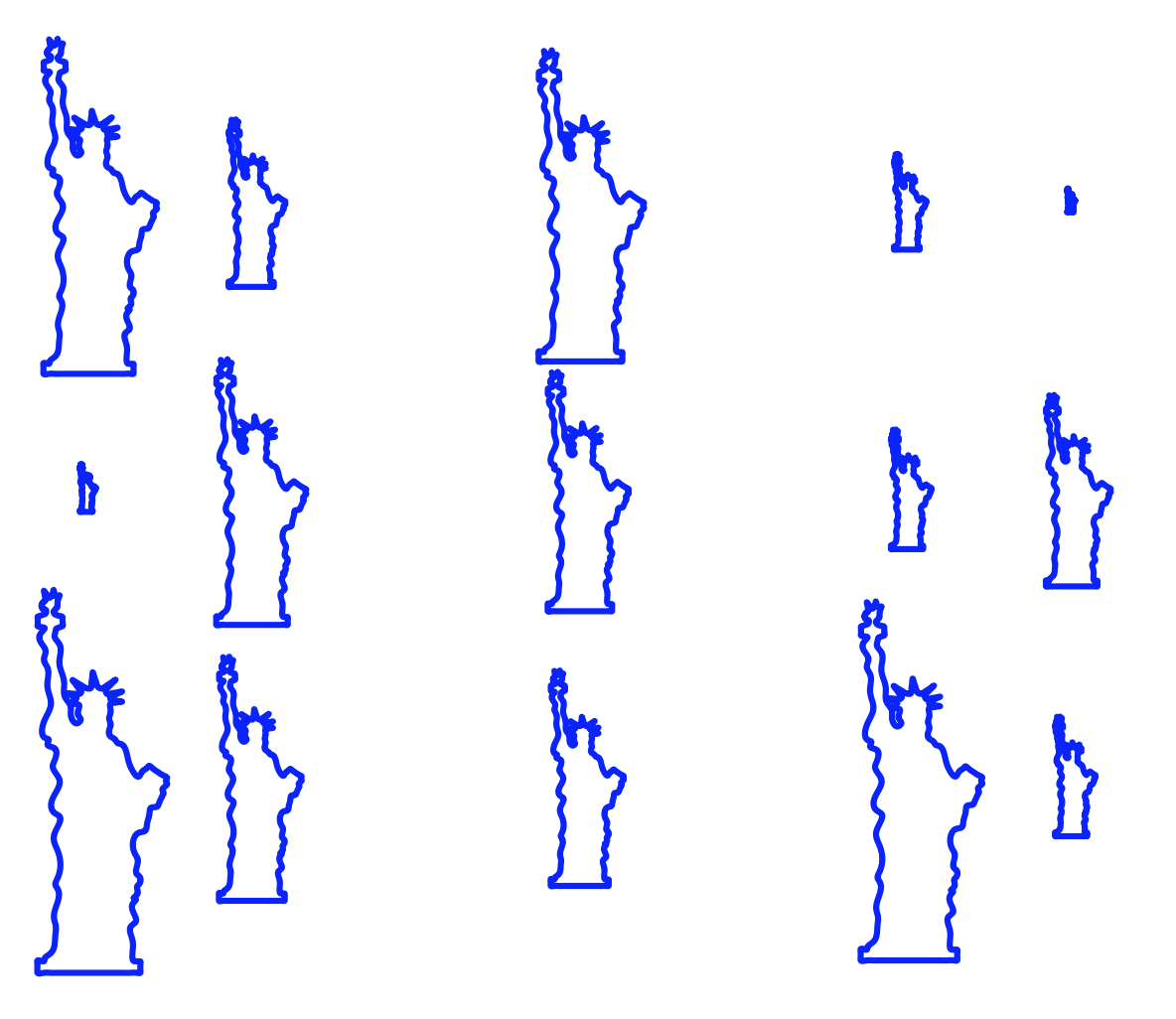} \end{center}& \begin{center} \begin{tabular}{l} \includegraphics[height=1.5cm]{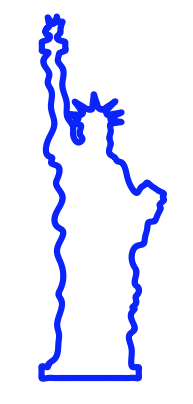}\\  length = 1 \end{tabular}\end{center} &\begin{center} \begin{tabular}{l} \includegraphics[height=2cm]{Liberty_bleue}\\ enclosed area =1 \end{tabular}\end{center}\\
\hline
\end{tabular}
\caption{ Examples of group actions on 2D simple closed curves and different choices of sections of the corresponding fiber bundle.}
 		\label{Tabular_group_action}	
\end{figure*}

\section{Introduction, motivation, and a simple example}

 In this paper we will describe three ways to think about geometry on a quotient space of a trivial principal bundle,
with application to shape space. The first is the standard approach via quotients by a group and a Riemannian submersion. The second is by considering
a particular global section of the bundle and inducing a metric by isometric immersion. The third is newer and consists of specifying a normal bundle
complementary to the vertical bundle, projecting the metric onto the normal bundle, and taking the quotient of the resulting degenerate metric; we refer to this as the gauge-invariant approach. We will begin with some motivations about our main concern of shape space before presenting an explicit example in finite dimensions to fix ideas. Then we describe the three basic methods as (I), (II), and (III), and finally we discuss how to get from one to another and the meaning of gauge invariance.

In order to explain the ideas of the present paper in a simple way, we will consider the contour of the Statue of Liberty appearing in different layouts in Fig.~\ref{Tabular_group_action}. Imagine a camera that scans a photo of the Statue and needs to recognize it regardless of how the photo is held (at any distance, position, or rotation angle). If we require the photo to be held perpendicular to the camera lens, then the transformations are rotations, translations, and rescalings. These are all shape-preserving  and we would like the scanner to be able to detect the shape independently of them.
We can handle this in two ways:
\begin{enumerate}
\item either one groups together the photos that are transformations of each other,
\item or one specifies a preferred choice of position in space and/or scale as a representative.
\end{enumerate}
The first option consists of considering the \textit{orbit} of the photo under the group; the second option consists of considering a preferred \textit{section} of the quotient space of curves modulo the group. These two different ways of thinking about shapes modulo a given group of transformations are illustrated in Fig.~\ref{Tabular_group_action}. In the second column we show some examples of the photos in the same orbit under the action of the group. In the third and fourth columns, a representative of this orbit is singled out.

Then, when considering contours of objects, another group acting by shape-preserving transformations is the group of reparameterizations of the contour. Making the analysis invariant by the group of reparameterizations
 is a much more difficult problem than that in the previous paragraph, but it has the same essential nature, and is the main source of motivation for us.

\textbf{Corresponding mathematical objects.}
The mathematical picture to start with is the following: the group of shape-preserving transformations $G$ is acting on the space of curves or surfaces $\mathcal{F}$ and the shape space $\mathcal{S}$ that retain just the informations that we need is the \textit{quotient space} $\mathcal{S}:=\mathcal{F}/G$. The map that sends a curve to its orbit under the group $G$ is called the \textit{canonical projection} and will be denoted be $p:\mathcal{F}\rightarrow \mathcal{S}$. The \textit{orbit} of a element $f\in \mathcal{F}$ will also be denoted by $[f]\in\mathcal{F}/\mathcal{S}$, in particular $p(f) = [f]$ for any $f\in\mathcal{F}$. The triple $(p, \mathcal{F}, \mathcal{S})$ is a particular example of fiber bundle attached to a smooth action of a group on a manifold.

When we specify which procedure we follow to choose a representant of each orbit, one is selecting a preferred \textit{section} of the fiber bundle $p:\mathcal{F}\rightarrow \mathcal{S}$. A global section of the fiber bundle $p:\mathcal{F}\rightarrow \mathcal{S}$ is a smooth application $s:\mathcal{S}\rightarrow \mathcal{F}$, such that $p\circ s([f]) = [f]$ for any $[f]\in\mathcal{S}$. There is one-to-one correspondance between the shape space $\mathcal{S}$ and the range of $s$. Defining a global section of $p:\mathcal{F}\rightarrow \mathcal{S}$ is in fact defining a way to choose a preferred element in the fiber $p^{-1}([f])$ over $[f]$. In the case of the group of reparameterizations, it consists of singling out a preferred parameterization of each oriented shape.

\textbf{Why do we care about the distinction?} Depending on the representation of shape space as a quotient space or as a preferred section, shape analysis may give different results. Very often, curves or surfaces are centered and scaled as a pre-processing step. However, the procedure to center or scale the shapes may influence further analysis. For instance, a Statue of Liberty whose contour has a fractal behaviour will appear very small if scaling variability is taken care of by fixing the length of the curve to $1$ and will seem visually very different to analogous statues with smooth boundaries.

\begin{example}\label{heisenberg}
We begin with the simplest nontrivial example of the three methods we have in mind for producing a metric on the quotient space
by a group action, given a metric on the full space. Here our full space will be the Heisenberg group $\mathcal{F}\cong\mathbb{R}^3$ with the left-invariant metric
\vspace{-0.2cm}
\begin{equation}\label{heisenbergmetric}
ds^2 = dx^2 + dy^2 + (dz-y\,dx)^2
\vspace{-0.3cm}
\end{equation}
on it, while the group action is vertical translation in the $z$-direction by a real number, generated by the flow of the vector field $\xi=\partial_z$.  Hence the group is $G=\mathbb{R}$ under addition, and the quotient space is $\mathcal{S} = \mathbb{R}^2$ with projection $p(x,y,z)=(x,y)$. We will denote by $\{e_1, e_2\}$ the canonical basis of $\mathbb{R}^2$. The metric~\eqref{heisenbergmetric} is invariant under this action since none of the components depend on $z$.\\
\textbf{(I)} At every point the field $\xi$ is vertical since it projects to zero. Horizontal vectors are those orthogonal to this in the metric \eqref{heisenbergmetric}, and the horizontal bundle
     is spanned 
     by the fields $h_1 = \partial_x + y \, \partial_z$  and $h_2 = \partial_y$.
 This basis is special since $Dp(h_1) = e_1$ and $Dp(h_2) = e_2$, so we get the usual basis on the quotient.
 If the inner product on $\{e_1,e_2\}$ comes from the inner product on $\{h_1,h_2\}$, the result is the Riemannian submersion quotient metric
\vspace{-0.3cm}
\begin{equation}\label{spatialsubmersion}
ds^2 = du^2+dv^2.
\vspace{-0.3cm}
\end{equation}
\textbf{(II)} The second way to get a natural metric on the quotient $\mathcal{S}$ is to embed it back into $\mathcal{F}$ by a section $s\colon \mathcal{S}\to \mathcal{F}$, so that $p\circ s$ is the identity. All such sections are given by the graph of a function $(x,y,z)=s(u,v) = (u,v,\psi(u,v))$ for some $\psi\colon \mathbb{R}^2\to\mathbb{R}$, a choice of a particular representative $z=\psi(u,v)$ in the equivalence class $\pi^{-1}(u,v)$. The image of $s$ is a submanifold of $\mathcal{F}$ which we denote by $\mathcal{A}$,  and it inherits the isometric immersion metric
\vspace{-0.3cm}
\begin{equation}\label{spatialimmersion}
ds^2 = du^2 + dv^2 + \big[\psi_v \, dv + (\psi_u - v)\,du\big]^2.
\vspace{-0.3cm}
\end{equation}
\textbf{(III)} The third way to get a metric is to declare that movement in the $z$ direction will be ``free,'' and only movement transverse to the vertical direction will have some cost. This corresponds to specifying a space of normal vectors along each fiber (arbitrary except that it is transverse to the tangent vectors $\partial_z$). Any normal bundle is generated by the span of vector fields of the form $n_1 = \partial_x + \varphi_1(x,y) \, \partial_z$ and $n_2=\partial_y + \varphi_2(x,y)\,\partial_z$ for some functions $\varphi_1,\varphi_2\colon \mathbb{R}^2 \to\mathbb{R}$ independent of $z$ to ensure $G$-invariance. Again this basis is special since $\pi_*(n_1)=e_1$ and $\pi_*(n_2)=e_2$. To measure movement only in the normal direction, we define $g_{GI}(U,V) = g_{\mathcal{F}}\big(p_N(U), p_N(V)\big)$ for any vectors $U$ and $V$, where $p_N$ is the projection onto the normal bundle parallel to the vertical direction. This results in the degenerate metric
\vspace{-0.3cm}
\begin{equation}\label{spatialgauge}
ds^2 = dx^2 + dy^2 + \big[\varphi_2 \, dy + (\varphi_1 - y)\,dx\big]^2.
\vspace{-0.3cm}
\end{equation}
This formula then induces a nondegenerate quotient metric on the quotient $\mathbb{R}^2$.

It is clear that (III) is the most general choice, and that both (I) and (II) are special cases. The metric \eqref{spatialgauge} matches \eqref{spatialsubmersion} when the normal bundle coincides with the horizontal bundle, and otherwise is strictly larger. Meanwhile the immersion metric (II) in \eqref{spatialimmersion} is strictly larger than (I), and no choice of $\psi$ will reproduce it since $\psi$ would have to satisfy $\psi_v = 0$ and $\psi_u = v$. This is a failure of integrability, see Section~\ref{Relationships}.
\end{example}

%
%
%
%
%

\section{Different methods to endow a quotient with a Riemannian metric}

In this section, we will suppose that we have at our disposal a Riemannian metric $g_{\mathcal{F}}$ on the space of curves or surfaces $\mathcal{F}$ we are interested in, and that this metric is invariant under a group of shape-preserving transformations $G$. In other words, $G$ preserves the metric $g_\mathcal{F}$, i.e. $G$ acts by isometries on $\mathcal{F}$. We will explain three different ways to endow the quotient space $\mathcal{S}:=\mathcal{F}/G$ with a Riemannian metric.

The action of a group $G$ on a space of curves or surfaces $\mathcal{F}$ will be denoted by a dot. For instance, if $G$ is the group of translations acting on curves in $\mathbb{R}^2$,  $g\cdot F = F + C$ where $C$ is the constant function given by the coordinates of the vector of translation defined by the translation $g$. When $G$ is the group of reparameterizations, then $g\cdot F := F \circ g^{-1}$.

\textbf{(I) Quotient Riemannian metric.} The first way to endow the quotient space $\mathcal{S}:=\mathcal{F}/G$ with a Riemannian metric is through the \textit{quotient Riemannian metric}. We recall the following classical Theorem of Riemannian geometry.
\begin{theorem}[Riemannian submersion Theorem]
Let $\mathcal{F}$ be a manifold endowed with a Riemannian metric $g_\mathcal{F}$, and $G$ a Lie group acting on $\mathcal{F}$ in such a way that $\mathcal{F}/G$ is a smooth manifold. Suppose $g_{\mathcal{F}}$ is $G$-invariant and $T_F\mathcal{F}$ splits into the direct sum of the tangent space to the fiber and its orthogonal complement, i.e.,
\begin{equation}\label{addition}
\begin{split}
g_\mathcal{F}(X, Y) = g_\mathcal{F}(g\cdot X, g\cdot Y), \forall X, Y \in T\mathcal{F}, \forall g \in G, \\
T_F\mathcal{F} = \operatorname{Ker}(dp)_F \oplus \operatorname{Ker}(dp)_F^\perp, \forall F \in \mathcal{F},
\end{split}
\end{equation}
then there exists a unique Riemannian metric $g_{1,\mathcal{S}}$ on the quotient space $\mathcal{S} = \mathcal{F}/G$ such that the canonical projection $p: \mathcal{F} \rightarrow S$ is a Riemannian submersion, i.e. such that
$
dp: \operatorname{Ker}(dp)^\perp \rightarrow T\mathcal{S}
$
is an isometry.
\end{theorem}
In this Theorem, the space $\operatorname{Hor}:= \operatorname{Ker}(dp)^\perp$ is called the \textit{horizontal space} because it is defined as the orthogonal with respect to $g_{\mathcal{F}}$ of the \textit{vertical space} $\textrm{Ver} := \operatorname{Ker}(dp)$ (traditionally the fibers of a fiber bundle are depicted vertically). Condition~\eqref{addition} is added in order to deal with the infinite-dimensional case where, for weak Riemannian metrics, this identity is not automatic.

One way to understand the Riemannian submersion Theorem is the following: first, in order to define a Riemannian metric on the quotient space, one looks for a subbundle of $T\mathcal{F}$ which is in bijection with $T\mathcal{S}$. Since the vertical space is killed by the projection,  the transverse space to the vertical space given by the orthogonal complement is a candidate.  The restriction of the Riemannian metric on it defines uniquely a Riemannian metric on the quotient.


\textbf{(II) Riemannian metric induced on a smooth section.}
Now suppose that we have chosen a preferred smooth section $s:\mathcal{S}\rightarrow \mathcal{F}$ of the fiber bundle $p: \mathcal{F}\rightarrow S = \mathcal{F}/G$, for instance the space of arc-length parameterized curves in the case where $G$ is the the group of orientation-preserving reparameterizations, or the space of centered curves when $G$ is the group of translations. The smoothness assumption means that the range of $s$ is a smooth manifold of $\mathcal{F}$, like the space of arc-length parameterized curves in the space of parameterized curves. We will denote it by $\mathcal{A}:=s(\mathcal{S})$. By construction, there is a isomorphism between $\mathcal{S}$ and $\mathcal{A}$ which one can use to endow the quotient space $\mathcal{S}$ with the induced Riemannian structure on $\mathcal{A}$ by $\mathcal{F}$.
\begin{theorem}[Riemannian immersion Theorem]
Given a smooth section $s: \mathcal{S}\rightarrow \mathcal{F}$, there exists a unique Riemannian metric $g_{\mathcal{A}}$ on $\mathcal{A}:= s(\mathcal{S})$ such that the inclusion $\iota: \mathcal{A} \hookrightarrow \mathcal{F}$ is an isometry. Using the isomorphism $s: \mathcal{S}\rightarrow \mathcal{A}$,
there exists a unique Riemannian metric $g_{2, \mathcal{S}}$ on $\mathcal{S}$ such that $s: \mathcal{S} \rightarrow \mathcal{F}$ is an isometry.
\end{theorem}

%
\textbf{(III) Gauge invariant metric}.
Here we suppose that we have a vector bundle $\operatorname{Nor}$ over $\mathcal{F}$ which is a $G$-invariant subbundle of  $T\mathcal{F}$ transverse to the vertical bundle $\textrm{Ver} := \operatorname{Ker}(dp)$.  Using any $G$-invariant metric $g_{\mathcal{F}}$ on $\mathcal{F}$, one can define a $G$-invariant metric $g_{GI}$ on $\mathcal{F}$ that is degenerate along the fiber of the projection $p: \mathcal{F}\rightarrow \mathcal{S}$. We will explain the meaning of ``gauge invariance'' later.
\begin{theorem}
Let $g_{\mathcal{F}}$ be a $G$-invariant metric on $\mathcal{F}$ and  $\operatorname{Nor}\subset T\mathcal{F}$ be a $G$-invariant subbundle of $T\mathcal{F}$ such that
\begin{equation}\label{dec}
T_F\mathcal{F} =  \operatorname{Ker}(dp)_F \oplus  \operatorname{Nor}_F, \forall F \in \mathcal{F}.
\end{equation}
There exists a unique metric $g_{GI}$ on $T\mathcal{F}$ which coincides with $g_\mathcal{F}$ on $\operatorname{Nor}$ and is
degenerate exactly along the vertical fibers of $p: \mathcal{F}\rightarrow \mathcal{S}$. It induces a  Riemannian metric $g_{3, \mathcal{S}}$ on shape space $\mathcal{S}$ such that $dp: \operatorname{Nor} \rightarrow T\mathcal{S}$ is an isometry.
\end{theorem}

Since we want the inner product to be the same in $g_{GI}$ as in $g_{\mathcal{F}}$ when the vectors 
are normal, and zero if either vector is vertical,
we define $g_{GI}$ by simply projecting an arbitrary vector onto the normal bundle:
\begin{equation}\label{GI}
g_{GI}(X,Y) = g_{\mathcal{F}}\big(p_{\text{Nor}}(X), p_{\text{Nor}}(Y)\big),
\end{equation}
where $p_N: T_F\mathcal{F}\rightarrow \text{Nor}$ is the projection onto the normal bundle parallel to the vertical space.
This is nondegenerate on the quotient since the projection onto the quotient is an isomorphism when restricted to the normal bundle.

\begin{remark}
In the case where $\operatorname{Nor} = \operatorname{Hor}$, the Riemannian metric $g_{3, \mathcal{S}}$ coincides with the quotient metric $g_{1, \mathcal{S}}$. Another choice of $G$-invariant complement to the vertical space will give another Riemannian metric on the quotient space.
\end{remark}

%
%

\begin{example}
The main example for shape space consists of the elastic metric first defined in \cite{elasticoriginal} on the space of planar curves $\mathcal{F} = \{F\colon [0,1]\to\mathbb{R}^2\}$
by the formula
\begin{multline}\label{elasticmetricdef}
g^{a,b}_{\mathcal{F}}(h_1,h_2) = \int_0^1 \big[ a(D_sh_1,\tangent)(D_sh_2,\tangent) + b(D_sh_1,\normal)(D_sh_2,\normal)\big] \, ds,\\
F\in \mathcal{F},\; h_i\in T_F\mathcal{F}, ds = \lVert F'(t)\rVert dt,\; D_sh(t) = \tfrac{\dot{h}(t)}{\lVert \dot{F}(t)\rVert}, \; \tangent = \tfrac{\dot{F}}{\lVert \dot{F}\rVert}, \; \normal = \tangent^{\perp}.
\end{multline}
See \cite{bauersurvey} for a recent survey of its properties. We will follow \cite{PrestonWhips} and \cite{TumPres} below.

Our group $G$ is the (orientation-preserving) reparameterizations of all these curves, since we only care about the image $F[0,1]$, and the shape space is the quotient $\mathcal{F}/G$. At any $F\in \mathcal{F}$ the vertical space is $\Ver = \{ m\tangent \, \vert \, m\colon [0,1]\to \mathbb{R}\}$.
A natural section $s\colon \mathcal{S}\to \mathcal{A}\subset \mathcal{F}$ comes from parameterizing all curves proportional to arc length. The tangent space  $T\mathcal{A}$ to the space of arc-length parameterized curves is the space of vector fields $w$ along $F$ such that $w'\cdot \tangent = 0$.
The horizontal bundle in the
metric \eqref{elasticmetricdef} is given at each $F$ by the space of vector fields $w$ along $F$ such that $\frac{d}{dt} (w'\cdot \tangent) - \tfrac{b}{a} \kappa(w'\cdot \normal)=0$ for all $t\in [0,1]$, where $\kappa$ is the curvature function.
Hence computing the projections requires solving an ODE.

A much simpler normal bundle is obtained by just taking the pointwise normal, i.e., using $\Nor :=  \{ \Phi \normal\, \vert \, \Phi\colon [0,1]\to\mathbb{R}\}$, where now the tangential and normal projections can be computed without solving an ODE.

Instead of parameterizing by arc length we can choose other special parameterizations to get the section; for example using speed proportional to the curvature of the shape as in \cite{Tum3}. An example of application to action recognition is given in \cite{Drira}. Similar metrics can be defined on surfaces in $\mathbb{R}^3$ to get two-dimensional shape spaces; see for example \cite{Tum1,Tum2}.
\end{example}

\begin{remark}
In the infinite-dimensional case, it is not always possible to find a complement to the vertical space $\operatorname{Ker}(dp)$ as in \eqref{dec}. An example of this phenomemon is provided by shape spaces of non-linear flags (see \cite{Ioana}). In this case, one has to work with the quotient vector spaces $T_{F}\mathcal{F}/\operatorname{Ker}(dp)_F$. See \cite{futurepaper} for this more general case.
\end{remark}
\section{Relationships of the $3$ methods and Gauge Invariance}\label{Relationships}

\subsection{Converting between (I), (II), and (III).}
We have seen in Example \ref{heisenberg} that in some cases the three metrics coincide when we start with the \emph{same}
base metric $g_{\mathcal{F}}$, but typically they do not. However if we allow the metric on $\mathcal{F}$ to change, we can convert any
metric of the form (I), (II), or (III) into a metric of the other forms. Here we demonstrate how to do it.

\textbf{(I)$\Rightarrow$(III).} If we start with a quotient Riemannian submersion metric arising from $g_{\mathcal{F}}$, how do we get a gauge-invariant metric? We simply define the normal bundle $\operatorname{Nor}$ to be the horizontal bundle $\operatorname{Hor}$ of vectors
orthogonal in $g_{\mathcal{F}}$ to the vertical bundle, and use the projection $p_{\operatorname{Nor}}$ as in \eqref{GI}.
The new metric $g_{GI}$ on $\mathcal{F}$ will be degenerate but will produce the same metric on the quotient.

\textbf{(II)$\Rightarrow$(III).} If we start with a section $s$ that embeds the quotient $\mathcal{S}$ into a submanifold $\mathcal{A}$ of $\mathcal{F}$, how do we obtain a gauge-invariant metric? Here we define the normal bundle $\operatorname{Nor}$ to be the tangent bundle of the $\mathcal{A}$ and proceed as in \eqref{GI}. Again the new degenerate metric on $\mathcal{F}$ will agree with the induced metric on $\mathcal{A}$ (and in particular be nondegenerate there).

\textbf{(I)$\Rightarrow$(II), (III)$\Rightarrow$(II)} As in Example \ref{heisenberg}, a given normal bundle (in particular a horizontal bundle from a metric) may not be the tangent bundle of any manifold due to failure of integrability; hence there may not be any way to express a particular instance of (I) or (III) as a version of (II) with the horizontal space equal to the tangent space of a particular section. The Frobenius integrability condition for the bundle, which is equivalent to the curvature of the corresponding connection vanishing, has to be satisfied in order to build a section tangent to the horizontal space. If one does not require the section to be horizontal, we may proceed as in \cite{TumPres} to pull-back the metric from the quotient to the particular section.

\textbf{(III)$\Rightarrow$(I), (II)$\Rightarrow$(I).} Both these cases come from the same metric $g_{\mathcal{F}}$, and as in Example \ref{heisenberg} they may not coincide. However if we are willing to consider a \emph{different} metric on $\mathcal{F}$, then we can go from method (III) to method (I) and still have the same metric on the quotient. Given a $G$-invariant normal bundle $\Nor$ and a $G$-invariant metric $g_{2,\mathcal{F}}$ which generates the degenerate metric  $g_{GI}$ on $\mathcal{F}$ and a metric $g_{\mathcal{S}}$ by method (III) via \eqref{GI}, we define a new metric $g_{1,\mathcal{F}}$ by choosing any $G$-invariant Riemannian metric on the tangent bundle $\operatorname{Ver}$ to the fiber and by declaring that $\operatorname{Ver}$ and $\operatorname{Nor}$ are orthogonal.
By construction, the subbundle $\Nor$ is the horizontal bundle in this metric, and we are in case (I). The same construction works to go from (II) to (I), using the above to get from (II) to (III).

\vspace{-0.5cm}
\subsection{The meaning of gauge invariance.}
In geometry on shape space, we are often interested primarily in finding minimizing paths between shapes, and a common algorithm is to construct an initial path between shapes and shorten it by some method.  Paths are typically easy to construct in $\mathcal{F}$ and difficult to construct directly in the quotient space $\mathcal{S}$. A section $s\colon \mathcal{S}\to \mathcal{F}$ makes this simpler, but especially if the image $\mathcal{A}=s(\mathcal{S})$ is not flat, it can be difficult to keep the shortening constrained on that submanifold. Our motivating example is when the shapes are parameterized by arc length as in \cite{TumPres}, with a path-straightening or gradient descent algorithm based on the metric $g_{\mathcal{F}}$: the optimal reduction gives intermediate curves that are typically no longer parameterized by arc length, and in fact the parameterizations can become degenerate. As such we may want to apply the group action of reparameterization independently on each of the intermediate shapes to avoid this breakdown.

If $G$ is a group acting on a space $\mathcal{F}$, define the \emph{gauge group} as the group $\mathcal{G} := \{g:[0,1]\rightarrow G\}$ of paths in $G$ acting on the space of paths $\gamma\colon [0,1]\to \mathcal{F}$ in $\mathcal{F}$ by the obvious formula
 $(g\cdot \gamma)(t) = g(t)\cdot\gamma(t)$, i.e., pointwise action in the $t$ parameter. We would like the metric on $\mathcal{F}$ to have the property that lengths of paths are invariant under this action, which essentially allows us to change the section $s$ ``on the fly'' if it's convenient. Since this involves pushing the path in the direction of the fibers, it is intuitively clear that the metric will need to be degenerate in those directions. We have the following proposition, whose proof we defer to \cite{futurepaper}.

\begin{proposition}
The length of a path $[\gamma]$ in $\mathcal{S}$ measured with the quotient metric $g_{1, \mathcal{S}}$ is equal to the length of any lift $\gamma$ of $[\gamma]$ in $\mathcal{F}$ measured with the gauge invariant metric $g_{GI}$.
\end{proposition}

\begin{remark}
The length of $\gamma$ measured with the metric $g_{\mathcal{F}}$ differs from the length of $[\gamma]$ unless $\gamma$ is also horizontal. In this case, $\gamma = g_0\cdot \gamma_0$ for a fixed $g_0\in G$.
\end{remark}

\vspace{-0.5cm}
\section{Conclusion}
\vspace{-0.2cm}

We have shown how to view the problem of constructing a Riemannian metric on a principal bundle quotient such as
shape space in three different ways. Any desired metric on the quotient can be viewed as any one of the three depending
on what is computationally convenient. The gauge-invariant method (III) is the most general and flexible, capturing the other
two more familiar methods as special cases; it has the advantage that it is convenient for shape space computations without
needing to work on the difficult shape space explicitly.
%
%
%

\begin{thebibliography}{8}


\bibitem{bauersurvey}
M. Bauer, N. Charon, E. Klassen, S. Kurtek, T. Needham, T. Pierron,
\textit{Elastic metrics on spaces of Euclidean curves: theory and algorithms}, 2022,
\textcolor{blue}{\href{https://arxiv.org/pdf/2209.09862.pdf}{arXiv:2209.09862v1}}

\bibitem{Ioana} I. Ciuclea, A.B. Tumpach, C. Vizman, \textit{Shape spaces of non-linear Flags}, submitted to GSI23.

\bibitem{Drira} H. Drira, A.B.Tumpach, M. Daoudi, \textit{Gauge invariant framework for trajectories analysis},
DIFFCV Workshop 2015,
DOI: 10.5244/C.29.DIFFCV.6, \textcolor{blue}{\href{http://math.univ-lille1.fr/~tumpach/Site/research_files/GIF_curve.pdf}{GIFcurve.pdf}}

\bibitem{elasticoriginal} W. Mio, A. Srivastava and S. H. Joshi, \textit{On shape of plane elastic curves}, International Journal
of Computer Vision, \textbf{73} (2007), 307--324.

\bibitem{PrestonWhips} S.C. Preston, \textit{The geometry of whips}, Ann. Global Anal. Geometry, \textbf{41}(3) (2012).

\bibitem{Tum1} A.B. Tumpach,\textit{ Gauge invariance of degenerate Riemannian metrics}, Not. Amer. Math. Soc., April 2016. \textcolor{blue}{\href{http://math.univ-lille1.fr/~tumpach/Site/research_files/Notices_full.pdf}{Notices\_full.pdf}}

\bibitem{Tum2} A.B. Tumpach, H. Drira, M. Daoudi, A. Srivastava, \textit{Gauge invariant framework for shape analysis of surfaces}.
IEEE Trans. Pattern Anal. Mach. Intell.,
\textbf{38}(1), (2016). 


\bibitem{Tum3} A.B. Tumpach, \textit{On canonical parameterizations of $2D$-shapes}, submitted to GSI23.

\bibitem{TumPres} A.B. Tumpach, S.C. Preston, \textit{Quotient elastic metrics on the manifold of arc-length parameterized plane curves}, J.  Geom. Mech., \textbf{9}(2) (2017), 227--256.

\bibitem{futurepaper} A.B. Tumpach, S.C. Preston, \textit{Riemannian metrics on shape spaces: comparison of different constructions}, in progress.

\end{thebibliography}
%

\end{document}